% To: Andrzej <shlhetal@micron.net>
% Subject: Forwarded mail....
% Date: Wed, 25 Oct 2000 13:32:02 +0200 (IST)
% From: Saharon Shelah <shelah@math.huji.ac.il>
% Mime-Version: 1.0
% Content-Description: 
% X-sliced-and-diced-by: 'savemail' 0.4, Jul 2000

\input amssym

\magnification=\magstep1
\hsize=15,4truecm
\vsize=23.5truecm
\mathsurround=1pt

\def\chapter#1{\par\bigbreak \centerline{\bf #1}\medskip}

\def\section#1{\par\bigbreak {\bf #1}\nobreak\enspace}

\def\sqr#1#2{{\vcenter{\hrule height.#2pt
      \hbox{\vrule width.#2pt height#1pt \kern#1pt
         \vrule width.#2pt}
       \hrule height.#2pt}}}

\def\k{\kappa}
\def\o{\omega}
\def\s{\sigma}
\def\d{\delta}

\def\l{\lambda}
\def\r{\rho}

\def\a{\alpha}
\def\b{\beta}

\def\t{\tau}
\def\g{\gamma}

\def\n{\eta}

%=====================

\def\A{{\cal A}}
\def\B{{\cal B}}

\def\P{{\cal P}}

%========================

%\def\pm{\buildrel{\scriptscriptstyle +}\over {\scriptscriptstyle -}}
\def\th #1 #2. #3\par\par{\medbreak{\bf#1 #2.
\enspace}{\sl#3\par}\par\medbreak}
\def\rem #1 #2. #3\par{\medbreak{\bf #1 #2.
\enspace}{#3}\par\medbreak}
\def\proof{{\bf Proof}.\enspace}
\def\sqr#1#2{{\vcenter{\hrule height.#2pt
      \hbox{\vrule width.#2pt height#1pt \kern#1pt
         \vrule width.#2pt}
       \hrule height.#2pt}}}
\def\eop{\mathchoice\sqr34\sqr34\sqr{2.1}3\sqr{1.5}3}

%====================================================================%
                                                                     %
% A macro for making references and blocks.                          %
                                                                     %
\newdimen\refindent\newdimen\plusindent                              %
\newdimen\refskip\newdimen\tempindent                                %
\newdimen\extraindent                                                %
                                                                     %
%\refskip has to be determined by the user! Otherwise \parindent is  %
%used, in accordance with \item.                                     %
                                                                     %
\def\ref#1 #2\par{\setbox0=\hbox{#1}\refindent=\wd0                  %
\plusindent=\refskip                                                 %
\extraindent=\refskip                                                %
\advance\extraindent by 30pt                                         %
\advance\plusindent by -\refindent\tempindent=\parindent %           %
\parindent=0pt\par\hangindent\extraindent %                          %
{#1\hskip\plusindent #2}\parindent=\tempindent}                      %
\refskip=\parindent                                                  %
                                                                     %
%====================================================================%

\def\empty{\emptyset}

\def\raj{\restriction}

\def\nda{\mathrel{\lower0pt\hbox to 3pt{\kern3pt$\not$\hss}\downarrow}}
\def\nbot{\mathrel{\lower0pt\hbox to 4pt{\kern3pt$\not$\hss}\bot}}
\def\ekom{\mathrel{\lower3pt\hbox to 0pt{\kern3pt$\sim$\hss}\mapsto}}

\def\T{\Theta}
\def\anR{\mathrel{\lower1pt\hbox to 2pt{\kern3pt$R$\hss}\not}}
\def\anoR{\mathrel{\lower1pt\hbox to 2pt{\kern3pt$\overline{R}$\hss}\not}}

\def\noR{\anoR\ \ }
\def\anRm{\mathrel{\lower1pt\hbox to 2pt{\kern3pt$R^{-}$\hss}\not}}
\def\nRm{\anRm\ \ \ \ }

\null
\vskip 2truecm
\centerline{\bf CONSTRUCTING STRONGLY EQUIVALENT NONISOMORPHIC}
\centerline{\bf MODELS FOR UNSUPERSTABLE THEORIES, PART A}
\vskip 1truecm
\centerline{Tapani Hyttinen and Saharon Shelah$^{*}$}
\vskip 3truecm

\chapter{Abstract}

We study how equivalent nonisomorphic
models an unsuperstable theory can have. We measure the equivalence
by Ehrenfeucht-Fraisse games. This paper continues the work
started in [HT].

\chapter{1. Introduction}

In [HT] we looked how equivalent nonisomorphic models first-order theories
can have i.e. we tried to strengthen S.Shelah's nonstructure
theorems. We used Ehrenfeucht-Fraisse games to measure the equivalence
(see Definition 2.2 below). If the theory is unstable or
it has OTOP or it is superstable with DOP then we were able to
prove maximal results by assuming strong cardinal assumptions.
We showed that if $\l^{<\l}=\l$ then
there is a model $\A$ of the theory
such that $\vert\A\vert =\l$ and
for all $\l^{+},\l$-trees $t$ there is a model
$\B$ such that $\vert\B\vert =\l$, $\A\not\cong\B$ and
$\exists$ has a winning strategy in the Ehrenfeucht-Fraisse
game $G^{2}_{t}(\A ,\B )$.

By assuming only that the theory is unsuperstable we were not
able to say much if we tried to measure the equivalence
by the length of
Ehrenfeucht-Fraisse games in which $\exists$ has a winning strategy.
But if instead, we measured
the equivalence by the length of
Ehrenfeucht-Fraisse games in which $\forall$
does not have a winning strategy,
then we were able to
get rather strong results.

In this paper we look the unsuperstable case again. We
measure the equivalence by the length of
Ehrenfeucht-Fraisse games in which $\exists$ has a winning strategy.
We study $\l^{+},\k +1$-trees (see Definition 2.1)
and give a rather complete answer to the question: how equivalent
nonisomorphic $\l^{+},\k +1$-trees can there be? In Chapter 3 we
show that if $\l =\mu^{+}$, $cf(\mu )=\mu$, $\k =cf(\k )\le\mu$
and $\l^{<\k}=\l$ then there are

\vskip 2truecm
\noindent
{\it * partially supported by the United States Israel binational
science foundation, publ. 474}
\vfill
\eject

\noindent
$\l^{+},\k +1$-trees $I_{0}$ and
$I_{1}$ such that $\vert I_{0}\vert\cup\vert I_{1}\vert\le\l^{\k}$,
$I_{0}\not\cong I_{1}$ and
$$I_{0}\equiv^{\l}_{\mu\times\k}I_{1}$$
(see Definition 2.2 and Definition 2.4 (iii)).
Instead of two such trees it is possible
to get $2^{\l}$ such trees.

In chapter 4 we show that if in addition
$\l\in I[\l ]$ then the result of Chapter 3 is best possible.

As in [HT], this implies that essentially the same is true also for the models
of the canonical example of unsuperstable theories.

In [HS] we will prove the results of chapter 3 for unsuperstable
theories in general.

This paper was born during the first author's visit to the second
author at Rutgers University. The first author wishes to express
his gratitude to Rutgers University for the hospitality shown to him
during the visit.

\chapter{2. Basic definitions}

In this chapter we define the basic concepts we shall use.

\th 2.1 Definition. Let  $\l$ be a cardinal and $\a$ an ordinal.
Let $t$ be a tree (i.e. for all $x\in t$, the set $\{ y\in t\vert\ y<x\}$
is well-ordered by the ordering of $t$).
If $x,y \in t$ and $\{ z \in t \mid z < x \} = \{ z \in t \mid z < y \}$,
then we denote $x \sim y$, and the equivalence class
of $x$ for $\sim$ we denote $[x]$.
By a $\l, \a$-tree $t$ we mean a
tree which satisfies:

(i) $\vert [x] \vert < \l$ for every $x \in t$;

(ii) there are no branches of length $\ge \a$ in $t$;

(iii) $t$ has a unique root;

(iv) if $x,y \in t$, $x$ and $y$ have no immediate predecessors
and $x\sim y$, then
$x=y$.

If $t$ satisfies only (i), (ii) and (iii) above, we say that $t$ is
a wide $\l ,\a$-tree.

Note that in a $\l,\a$-tree each ascending sequence of a limit length
has at most one supremum, but in a wide $\l ,\a$-tree an ascending sequence
may have more than one supremum.

\th 2.2 Definition. Let $t$ be a tree and $\k$ a cardinal.
The  Ehrenfeucht-Fraisse game
of length $t$ between models $\A$ and $\B$,
$G^{\k}_{t}(\A, \B)$, is the following.
At each move $\a$:

(i) player $\forall$ chooses $x_\a \in t$, $\k_{\a}<\k$ and
either $a_\a^\b \in \A$, $\b <\k_{\a}$ or $b_\a^\b \in \B$,
$\b <\k_\a$, we will denote this sequence of elements of $\A$ or
$\B$
by $X_{\a}$;

(ii) if $\forall$ chose from $\A$ then
$\exists$ chooses $b_\a^\b \in \B$, $\b <\k_\a$, else
$\exists$ chooses
$a_\a^\b \in \A$, $\b <\k_\a$, we will denote this sequence by $Y_{\a}$.
\medskip
\noindent
$\forall$ must move so that $(x_\b)_{\b \le \a}$
form a strictly increasing sequence in $t$.
$\exists$ must move so that
$\{ (a_\g^\b, b_\g^\b) \vert \g \le \a , \b <\k_\g \}$
is a partial isomorphism from $\A$ to $\B$.
The player who first has to break the rules loses.

We write $\A\equiv^{\k}_{t}\B$ if $\exists$ has a winning strategy
for $G^{\k}_{t}(\A ,\B )$.

\th 2.3 Remark. Notice that the Ehrenfeucht-Fraisse game
$G^{\k}_{t}(\A ,\B )$ need not be determined, i.e. it may
happen that neither $\exists$ nor $\forall$ has a winning
strategy for $G^{\k}_{t}(\A ,\B )$ (see [MSV]).

\th 2.4 Definition. Let $t$ and $t'$ be trees.

(i) If $x \in t$, then $pred(x)$ denotes the sequence $(x_\a)_{\a < \b}$
of the predecessors of $x$, excluding $x$ itself,
ordered by $<$. Alternatively, we consider $pred(x)$ as a set.
The notation $succ(x)$ denotes the set of immediate
successors of $x$.
If $x,y \in t$ and there is $z$,
such that $x,y \in succ(z)$, then we say that
$x$ and $y$ are brothers.

(ii) By $t^{<\a}$ we mean the set
$$\{ x\in t\vert\ \hbox{\sl the order type of}\ pred(x)
\ \hbox{\sl is}\ <\a\} .$$
Similarly we define $t^{\le\a}$.

(iii) If $\a$ and $\b$ are ordinals then by
$\a +\b$ and $\a\times\b$ we mean ordinal sum and product (see [Je]).
Notice that ordinals are also trees.

\chapter{3. On nonstructure of trees of fixed height}

In this chapter we will assume that $\l =\mu^{+}$,
$cf(\mu )=\mu$, $\k =cf(\k )\le\mu$ and $\l^{<\k}=\l$.

Let $I^{+}_{n}=\{ \n\in\ ^{\le\k}\l\vert\ \n (0)=n\} -\{ ()\}$
and $I^{-}_{n}=\{ \n\in\ ^{<\k}\l\vert\ \n (0)=n\} -\{ ()\}$,
$n=0,1$. We consider these as trees ordered by initial segment
relation. Because for all $\d\le\k$, $(I^{+}_{n})^{<\d}=(I^{-}_{n})^{<\d}$
(see Definition 2.4), we denote this set by $I^{<\d}_{n}$ and
similarly we define $I^{\le\d}_{n}=(I^{+}_{n})^{\le\d}$
for all $\d <\k$.

If $\n\in I^{+}_{0}$ and $\xi\in I^{+}_{1}$
then we write $\n R^{-}\xi$ and $\xi R^{-}\n$ iff
$\n (j)=\xi (j)$ for all
$0<j<min\{ length(\n ),length(\xi )\}$ even.
For all $i<\k$ odd, we define $P_{i}$ to be the set of all
$\n\in I^{-}_{0}$ such that
$length(\n )=i$. Let $P=\bigcup \{P_{i}\vert\ i<\k ,
\ i\ \hbox{\rm odd}\}$

\th 3.1 Lemma. There is a partition
$\{ S_{\n}\vert\ \n\in P\}$ of $\l$ such
that for all $\n\in P$

(i) $\{ \d\in S_{\n}\vert\ cf(\d )=\mu\}$ is stationary;

(ii) if $\d\in S_{\n}$ and $cf(\d )=\mu$ then
$\d =sup(\d\cap S_{\n})$.

\proof Because $\vert P\vert =\l$ we can find a partition
of $\{ \a <\l\vert\ cf(\a )=\mu\}$ which satisfies (i). Let this
partition be $\{ S'_{\n_{\g}}\vert\ \g <\l\}$, where
$\{\n_{\g}\vert\ \g<\l\}$ is an enumeration of $P$.
Let $\{\s_{\g}\vert\ \g<\l\}$ be an enumeration of
$\{ \a <\l\vert\ cf(\a )=\mu\}$ so that
if $\s_{\g}>\s_{\g'}$ then $\g >\g'$.
We may assume that if $\d\in S'_{\n_{\g}}$,
$\g\ne 0$, then
$\d >\s_{\g}$. By induction on $\a\le\l$ we define sets $S^{\a}_{\n_{\g}}$.
Let $S^{0}_{\n_{0}}=S'_{\n_{0}}\cup\s_{0}$ and for all $\g >0$,
$S^{0}_{\n_{\g}}=S'_{\n_{\g}}$.
If $\a$ is limit ordinal and $cf(\a )\ge\mu$,
then we define $S^{\a}_{\n_{\g}}=
\bigcup_{\b <\a}S^{\b}_{\n_{\g}}$ for all $\g <\l$.
Assume $\a$ is successor or limit ordinal with $cf(\a )<\mu$.
Let $\s'_{\a}=\cup_{\d <\a}\s_{\d}$.
Then we choose $S^{\a}_{\n_{\g}}$ so that (a)-(f)
below are satisfied:

(a) $\bigcup_{\d <\a}S^{\d}_{\n_{\g}}\subseteq S^{\a}_{\n_{\g}}$,

(b) $S^{\a}_{\n_{\g}}\cap S^{\a}_{\n_{\g'}}=\empty$ if
$\g\ne\g'$,

(c) $\s_{\a}\subseteq\bigcup_{\g <\l}S^{\a}_{\n_{\g}}$,

(d) $S^{\a}_{\n_{\g}}-\s_{\a}=S^{0}_{\n_{\g}}-\s_{\a}$ for all
$\g<\l$,

(e) if $\s_{\a}\in S'_{\n_{\g}}$ then
$\s_{\a}=sup(\s_{\a}\cap S^{\a}_{\n_{\g}})$,

(f) if $\g\le\a$ then $(\s_{\a}-\s'_{\a})\cap S^{\a}_{\n_{\g}}\ne\empty$.

Then clearly $S_{\n_{\g}}=S^{\l}_{\n_{\g}}$, $\g<\l$, is a partition
of $\l$ and (i) is satisfied. We show that also (ii) is satisfied:
If $\s_{\d}\in S_{\n_{\g}}$ and $\d$ is successor or limit with
$cf(\d )<\mu$ then by (e) $\s_{\d}=sup(\s_{\d}\cap S_{\n_{\g}})$.
Otherwise we know that $\s_{\d}>\s_{\g}$ i.e. $\d >\g$ and
$sup\ \{\s_{\b}\vert\ \b <\d\} =\s_{\d}$. By (f) this implies that
$\s_{\d}=sup(\s_{\d}\cap S_{\n_{\g}})$.
$\eop$

\th 3.2 Definition. We define a relation $R\subseteq (I_{0}^{+}-I_{0}^{-})
\times (I_{1}^{+}-I_{1}^{-})$.
Let $\n\in I_{0}^{+}-I_{0}^{-}$ and
$\xi\in I_{1}^{+}-I_{1}^{-}$. Then $(\n ,\xi )\in R$ iff

(i) $\n R^{-}\xi$;

(ii) for every $j<\k$ odd, $\n$ and $\xi$ satisfy the following:
for all $\r\in P$, $\n (j)\in S_{\r}$ iff
$\xi (j)\in S_{\r}$ and
if
$\n (j)\not\in S_{\n\raj j}$,
then $\n (j)=\xi (j)$;

(iii) the set $W^{\k}_{\n ,\xi}$ is bounded in $\k$, where
$W^{\k}_{\n ,\xi}$ is defined in the following way:
Let $\d\le\k$,
$\n\in I^{+}_{0}-I^{<\d}_{0}$ and $\xi\in I^{+}_{1}-I^{<\d}_{1}$ then
$$W^{\d}_{\n ,\xi}=
\{ j<\d\vert\ j\ \hbox{\sl odd and}
\ \n (j)\in
S_{\n\raj j}\ \hbox{\sl and}$$
$$cf(\n (j))=
\mu\ \hbox{\sl and}
\ \xi (j)\ge\n (j)\} .$$

In order to simplify the notation we write $\n R\xi$ and
$\xi R\n$ for $(\n ,\xi)\in R$. Notice that by this we do not
try to claim that the relation is symmetric, in fact
it is antisymmetric, if $(\n ,\xi)\in R$ then always
$\n\in I_{0}^{+}-I_{0}^{-}$ and
$\xi\in I_{1}^{+}-I_{1}^{-}$. We also take liberty to write
$W^{\d}_{\xi ,\n}$ for $W^{\d}_{\n ,\xi}$ when it is convinient.

Our first goal in this chapter
is to prove the following theorem.
We will prove it in a sequence
of lemmas.

\th 3.3 Theorem. If $I_{0}$ and $I_{1}$ are such that

(i) $I^{-}_{n}\subseteq I_{n}\subseteq I^{+}_{n}$, $n=0,1$

\noindent
and

(ii) if $\n R\xi$, $\n\in I^{+}_{0}$ and $\xi\in I^{+}_{1}$ then
$\n\in I_{0}$ iff $\xi\in I_{1}$,

\noindent
then $I_{0}\equiv^{\l}_{\mu\times\k}I_{1}$.

\relax From now on
in this chapter we assume
that $I_{0}$ and $I_{1}$ satisfy (i) and (ii) above.

\th 3.4 Definition. Let $\a <\k$.

(i) $G_{\a}$ is the family of all
partial functions $f$ satisfying:

(a) $f$ is a partial isomorphism from $I_{0}$ to $I_{1}$;

(b) $dom(f)$ and $rng(f)$ are closed under initial segments and
for some $\b <\l$ they are
included in $\{\n\in I_{0}^{+}\vert
\ \hbox{\sl for all}\ j<\k ,\ \n (j)<\b\}$ and
$\{\xi\in I_{1}^{+}\vert
\ \hbox{\sl for all}\ j<\k ,\ \xi (j)<\b\}$, respectively;

(c) if $f(\n )=\xi$ then $\n R^{-}\xi$;

(d) if $\n\in I_{0}$, $\xi\in I_{1}$, $f(\n )=\xi$ and
$length(\n )=j+1$, $j$ odd, then $\n$ and $\xi$ satisfy the following:
for all $\r\in P$, $\n (i)\in S_{\r}$ iff
$\xi (i)\in S_{\r}$ and if
$\n (j)\not\in S_{\n\raj j}$,
then $\n (j)=\xi (j)$;

(e) assume $\n\in I^{+}_{0}-I^{<\d}_{0}$ and
$\{\n\raj\g\vert\ \g <\d\}\subseteq dom(f)$ and
let $\xi =\bigcup_{\g <\d}f(\n\raj\g )$, then
$W^{\d}_{\n ,\xi}$
has order type $\le\a$;

(f) if $\n\in dom(f)$ then $\{\g <\l\vert\ \n\frown (\g)\in dom(f)\}
=\{\g <\l\vert\ f(\n )\frown (\g)\in rng(f)\}$ is an ordinal.

(ii) We define $F_{\a}\subseteq G_{\a}$ by replacing (f) above by

(f') if $\n\in dom(f)$ then $\{\g <\l\vert\ \n\frown (\g)\in dom(f)\}
=\{\g <\l\vert\ f(\n )\frown (\g)\in rng(f)\}$ is an ordinal
of cofinality $<\mu$.

\th 3.5 Definition. For $f,g\in G_{\a}$ we write $f\le g$ if
$f\subseteq g$ and if $\g <\d \le\k$,
$\n\in I^{+}_{0}-I^{<\d}_{0}$, $\n\raj\g\in dom(f)$,
$\n\raj (\g+1)\not\in dom(f)$,
$\n\raj j\in dom(g)$ for all $j<\d$
and $\xi =\bigcup_{j<\d}g(\n\raj j)$,
then $W^{\g}_{\n ,\xi}=W^{\d}_{\n ,\xi}$.

Notice that $f\le g$ is a transitive relation.

\th 3.6 Remark. Let $f\in G_{\a}$.

(i) We define $\overline{f}$ by
$$dom(\overline{f})=dom(f)\cup\{\n\in I_{0}\vert\ \n\raj\g\in
dom(f)\ \hbox{\sl for all}\ \g <length(\n )$$
$$\hbox{\sl and}
\ length(\n )\ \hbox{\sl is limit}\}$$
and if
$\n\in dom(\overline{f})-dom(f)$ then
$$\overline{f}(\n )=\bigcup_{\g <length(\n )}f(\n\raj\g).$$

(ii) If $f\in F_{\a}$ then $\overline{f}\in F_{\a}$ and
if $f\in G_{\a}$ then
$\overline{f}\in G_{\a}$.

\th 3.7 Lemma. Assume $\a <\k$, $\d\le\mu$, $f_{i}\in F_{\a}$ for all
$i<\d$ and $f_{i}\le f_{j}$ for all $i<j<\d$.

(i) $\bigcup_{i<\d}f_{i}\in G_{\a}$.

(ii) If $\d <\mu$ then
$\bigcup_{i<\d}f_{i}\in F_{\a}$ and $f_{j}\le\bigcup_{i<\d}f_{i}$
for all $j<\d$.

\proof Follows immediately from the definitions. $\eop$

\th 3.8 Lemma. If $\d<\k$, $f_{i}\in G_{i}$ for all $i<\d$ and
$f_{i}\subseteq f_{j}$ for all $i<j<\d$ then
$$\bigcup_{i<\d}f_{i}\in G_{\d}.$$

\proof Follows immediately from the definitions. $\eop$

\th 3.9 Lemma. If $f\in F_{\a}$ and $A\subseteq I_{0}\cup I_{1}$,
$\vert A\vert <\l$, then there is $g\in F_{\a}$
such that $f\le g$ and
$A\subseteq dom(g)\cup rng(g)$.

\proof Let $\n\in dom(f)$ and let
$$\{ i<\l\vert\ \n\frown (i)\in dom(f)\} =
\{ i<\l\vert\ f(\n )\frown (i)\in rng(f)\} =\d ,$$
$cf(\d )<\mu$, and let $\b >\d$.
We show first that there are $f^{\n\b}\in F_{\a}$
and $\g\ge\b$ such that $f^{\n\b}\ge f$, $cf(\g )<\mu$ and
$$\{ i<\l\vert\ \n\frown (i)\in dom(f^{\n\b})\} =
\{ i<\l\vert\ f(\n )\frown (i)\in rng(f^{\n\b})\} =\g .$$

Let $length(\n )=j$. If $j$ is even it is trivial to find $f^{\n\b}$ and
$\g$. So we assume that $j$ is odd. We choose $\g\ge\b$ so that
$cf(\g )<\mu$. For any $i\in\g -\d$ satisfying:

(i) $cf(i)=\mu$

\noindent
and

(ii) $i\in S_{\n}$,

\noindent
we choose $j_{i}\in i-\d$ so that $j_{i}\in S_{\n}$,
$cf(j_{i})<\mu$ and
if $i\ne i'$ then $j_{i}\ne j_{i'}$. These $j_{i}$ exist
because $sup\ i\cap S_{\n}=i$ and $i\ne\d$.

Then we define $f^{\n\b}(\n\frown (i))=f(\n )\frown (j_{i})$ and
$f^{\n\b}(\n\frown (j_{i}))=f(\n )\frown (i)$. For all other
$i\in\g -\d$ we let $f^{\n\b}(\n\frown (i))=f(\n )\frown (i)$.
It is easy to see that $f^{\n\b}\in F_{\a}$ and $f^{\n\b}\ge f$.

It is easy to see that we can choose $\n_{i}\in I_{0}$ and
$\b_{i}<\l$, $i<\mu$, so that the following functions are well-defined:

(i) $g_{o}=f$;

(ii) $g_{i+1}=(g_{i})^{\n_{i}\b_{i}}$;

(iii) $g_{i}=\overline{(\bigcup_{j<i}g_{j})}$, if $i$ is limit;

\noindent
and $A\subseteq dom(\bigcup_{i<\mu}g_{i})\cup rng(\bigcup_{i<\mu}g_{i})$.
Furthermore we can choose $\n_{i}$ and $\b_{i}$ so that
if $i\ne i'$ then $\n_{i}\ne\n_{i'}$. Then $g=\bigcup_{i<\mu}g_{i}$
is as wanted.
$\eop$

\th 3.10 Lemma. If $f\in G_{\a}$,
then there is $g\in F_{\a +1}$
such that $f\subseteq g$.

\proof Essentially as the proof of Lemma 3.9. $\eop$

Theorem 3.3 follows now easily from the lemmas above.

In the rest of this chapter we prove that there are trees
$I_{0}$ and $I_{1}$ which satisfy the assumptions of Theorem 3.3
and are not isomorphic. For this we use the following
Black Box. We define $H_{<\k^{+}}(\l )$ to be the smallest set H such
that

(i) $\l\subseteq H$

\noindent
and

(ii) if $x\subseteq H$ and $\vert x\vert\le\k$ then $x\in H$.

\th 3.11 Theorem. ([Sh3] Lemma 6.5) There is
$W=\{ (\overline{M}^{\a},\n^{\a})\vert\ \a <\a (*)\}$
such that:

(i) $\overline{M}^{\a}=(M^{\a}_{i}\vert\ i\le\k )$ is an increasing continuous
elementary chain of models belonging to $H_{<\k^{+}}(\l )$ and
$\n^{\a}\in\ ^{\k}\l$ is increasing;

(ii) $M^{\a}_{i}\cap\k^{+}$ is an ordinal, $\k +1\subseteq M^{\a}_{i}$,
$M^{\a}_{i}\in H_{<\k^{+}}(\n^{\a}(i))$,
$(M^{\a}_{j}\vert\ j\le i)\in M^{\a}_{i+1}$ and
$\n^{\a}\raj i\in M^{\a}_{i+1}$;

(iii) In the following game, $G(\k ,\l ,W)$, player $\forall$ does not have
winning strategy: The play lasts $\k$ moves, in the i-th move
$\forall$ chooses a model $M_{i}\in H_{<\k^{+}}(\l )$ and then
$\exists$ chooses $\g_{i}<\l$. $\forall$ must choose models
$M_{i}$, $i<\k$, so that
$(M_{i}\vert\ i\le\k )$ is an increasing continuous
elementary chain of models, $M_{i}\cap\k^{+}$ is an ordinal,
$\k +1\subseteq M_{i}$ and
$(M_{j}\vert\ j\le i)\in M_{i+1}$.
In the end $\exists$ wins the play
if for some $\a <\a (*)$, $\n^{\a}=(\g_{i}\vert\ i<\k )$ and
$M_{i}=M^{\a}_{i}$ for all $i<\k$;

(iv) $\n^{\a}\ne\n^{\b}$ for $\a\ne\b$.

Notice that in the game above $\forall$ can choose the similarity type
of models freely as long as other requirements are satisfied.

We define $I_{0}$ and $I_{1}$ with help of
$W$. We do this by defining $J_{\a}$, $\neg J_{\a}$, $K_{\a}$
and $\neg K_{\a}$ by induction
on $\a <\a (*)$ so that $J_{\a}\cap\neg J_{\a}=\empty$ and
$K_{\a}\cap\neg K_{\a}=\empty$
and then letting $I_{0}=I^{-}_{0}\cup\bigcup_{\a <\a (*)}J_{\a}$
and $I_{1}=I^{-}_{1}\cup\bigcup_{\a <\a (*)}K_{\a}$.
We assume that we have well-ordered $I^{+}_{0}-I^{-}_{0}$.

We say that $\a <\a(*)$ is active, if there is
$\n\in I^{+}_{0}-I^{-}_{0}$
such that $\a$ and $\n$ satisfy (i)-(vii) or
(i)-(v), (vi') and (vii') below.

(i) For all $i\le\k$,
the similarity type of $M^{\a}_{i}$ is $\{ \in ,I^{-}_{0},I^{-}_{1},g\}$
where $\in$ and $g$ are two-ary relation symbols and
$I^{-}_{0}$ and $I^{-}_{1}$ are unary relation symbols;

(ii) for all $i\le\k$,
$$M^{\a}_{i}\raj\{ \in ,I^{-}_{0},I^{-}_{1}\}\prec (H_{<\k^{+}}(\l ),\in ,
I^{-}_{0},I^{-}_{1});$$

(iii) for all $i<\k$, $\n\raj i\in M^{\a}_{i+1}$;

(iv) for all $i\le\k$, $M^{\a}_{i}\models$ "$g$ is an isomorphism from
$I^{-}_{0}$ to $I^{-}_{1}$";

(v) for all $\o\le i<\k$, if $i=\g +2k$ for some $\g$ limit
and $k<\o$ then $\n (i)=\n^{\a}(\g +k)$, and for all
$i<\o$, if $i=2k+2$ then $\n (i)=\n^{\a}(k)$;

let
$$\xi =\bigcup_{i<\k}g_{\a}(\n\raj i),$$
where $g_{\a}$ is the interpretation of $g$ in $M^{\a}_{\k}$,

(vi) $\n R^{-}\xi$

(vi') $\n\nRm \xi$

(vii) for all $i<\k$ odd, $\n(i)$ satisfies:

(a) $cf(\n(i))=\mu$ and
$\n (i)\in S_{\n\raj i}$;

(b) $M^{\a}_{\k}\models$ "the set
$\{\n\raj i\frown (j)\vert\ j<\n (i)\}\cup\{ g(\n\raj i)\frown (j)
\vert\ j<\n (i)\}$
is closed under $g$ and $g^{-1}$"

(vii') there is $j_{\n}<\k$ such that for all $i>j_{\n}$
odd the following holds:

(a) if $i=\g+4n+1$ for some limit ordinal $\g$ and $n\in\o$ then
$\xi (i)\in S_{\n\raj i}$

(b) if $i=\g+4n+3$ for some limit ordinal $\g$ and $n\in\o$ then
$\n (i)\in S_{\n\raj i}$, $cf(\n (i))=\mu$ and $\xi (i)\ge \n (i)$.

If $\a$ is active and there exists such $\n$ that
$\a$ and $\n$ satisfy (i)-(vii) above,
then we define $\n_{\a}$ to be
the least such $\n\in I^{+}_{0}-I^{-}_{0}$ in the well-ordering
of $I^{+}_{0}-I^{-}_{0}$.
Otherwise we let
$\n_{\a}$ to be the least $\n\in I^{+}_{0}-I^{-}_{0}$ in the well-ordering
of $I^{+}_{0}-I^{-}_{0}$ such that
$\a$ and $\n$ satisfy (i)-(v), (vi') and (vii') above.
Let
$$\xi_{\a} =\bigcup_{i<\k}g_{\a}(\n_{\a}\raj i),$$
where $g_{\a}$ is the interpretation of $g$ in $M^{\a}_{\k}$.
If $\a$ is active and $\n_{\a}\nRm \xi_{\a}$ then let
$j_{\a}=j_{\n_{\a}}$.

Let $\overline{R}$ be the transitive and reflexive closure of $R$.

\th 3.12 Lemma. If $\g$ is active
then $\n_{\g}\noR\xi_{\g}$.

\proof Clearly we may assume that $\n_{\g}R^{-}\xi_{\g}$.
For a contradiction assume, that there are
$\r_{0},...,\r_{n}$ such that $\r_{0}=\n_{\g}$,
$\r_{n}=\xi_{\g}$, for all $m<n$,
$\r_{m}R\r_{m+1}$ and for all $k<m\le n$, $\r_{k}\ne\r_{m}$.
We choose $i<\k$ so that

($\a$) $i$ is odd;

($\b$) for all $k<m\le n$, $\r_{k}\raj i\ne\r_{m}\raj i$;

($\g$) for all $m<n$, $W^{\k}_{\r_{m},\r_{m+1}}\subseteq i$.

\noindent
Because $\n_{\g} (i)\in S_{\n_{\g}\raj i}$ and $cf(\n_{\g} (i))=\mu$,
$\r_{1}(i)<\r_{0}(i)$ and $\r_{1}(i)\in S_{\n_{\g}\raj i}$.
By the definition of $R$, $\r_{2}(i)\in S_{\n_{\g}\raj i}$.
By ($\b$) above $\r_{2}(i)=\r_{1}(i)$ and $\r_{3}(i)=\r_{2}(i)$.
We can continue this and get $\r_{n}(i)=...=\r_{1}(i)$.
So $\n_{\g} (i)>g(\n_{\g}\raj (i+1))(i)$ which contradicts with
(vii)(b) in the definition of active. $\eop$

\th 3.13 Lemma. Let $\a$ and $\b$ be active, $\a\ne\b$,
$\xi_{\a}\overline{R}\xi_{\b}$
and $\n_{\a}\nRm \xi_{\a}$ then $\n_{\b}R^{-}\xi_{\b}$.

\proof For a contradiction assume $\n_{\b}\nRm \xi_{\b}$. By (vii') (a)
in the definition of active
we can find $i<\k$ odd such that $\xi_{\a}(i)\in S_{\n_{\a}\raj i}$  and
$\xi_{\b}(i)\in S_{\n_{\b}\raj i}$. By Definition 3.2 (ii) this implies
$\xi_{\a}\noR\xi_{\b}$, a contradiction.
$\eop$

\th 3.14 Lemma. Let $\a$ and $\b$ be active.

(i) If $\a\ne\b$ then $\n_{\a}\noR\n_{\b}$.

(ii) If $\n_{\a}R^{-}\xi_{\a}$
then for all active $\g$, $\n_{\a}\noR\xi_{\g}$.

\proof (i) By (vii) (a) and (or) (vii') (b) in the definition of active
there is $i<\k$ odd
such that $\n_{\a}(i)\in S_{\n_{\a}\raj i}$,
$\n_{\b}(i)\in S_{\n_{\b}\raj i}$ and $\n_{\a}\raj i\ne\n_{\b}\raj i$.
By Definition 3.2 (ii) this implies
$\n_{\a}\noR\n_{\b}$.

(ii) If $\g =\a$ the claim follows immediately from Lemma 3.12.
So assume $\g\ne\a$. We may also assume $\n_{\a}R^{-}\xi_{\g}$,
because otherwise we have proved the claim. Then
$\n_{\g}\nRm\xi_{\g}$. By (vii) (a) and (vii') (a)
in the definition of active we can find
$i<\k$ odd such that $\n_{\a}(i)\in S_{\n_{\a}\raj i}$,
$\xi_{\g}(i)\in S_{\n_{\g}\raj i}$ and $\n_{\a}\raj i\ne\n_{\g}\raj i$.
As above this implies
$\n_{\a}\noR\xi_{\g}$.
$\eop$

\th 3.15 Lemma. Let $\a$ and $\b$ be active. If $\n_{\a}\overline{R}\xi_{\b}$
then there is $l_{\a\b}<\k$ such that for all $i>l_{\a\b}$,
$i=\g +4k+3$, $\g$ limit and $k\in\o$, $\n_{\a}(i)>\xi_{\b}(i)$.

\proof By Lemma 3.14 (ii) we may assume $\n_{\a}\nRm\xi_{\a}$.
For a contradiction assume, that there are
$\r_{0},...,\r_{n}$ such that $\r_{0}=\n_{\a}$,
$\r_{n}=\xi_{\b}$, for all $m<n$,
$\r_{m}R\r_{m+1}$ and for all $k<m\le n$, $\r_{k}\ne\r_{m}$.
We choose $l_{\a\b}<\k$ so that

($\a$) $j_{\n_{\a}}<l_{\a\b}$;

($\b$) for all $k<m\le n$, $\r_{k}\raj i\ne\r_{m}\raj i$;

($\g$) for all $m<n$, $W^{\k}_{\r_{m},\r_{m+1}}\subseteq i$.

\noindent
Let $i>l_{\a\b}$,
$i=\g +4k+3$, $\g$ limit and $k\in\o$.
Because $\n_{\a} (i)\in S_{\n_{\a}\raj i}$ and $cf(\n_{\a} (i))=\mu$,
$\r_{1}(i)<\r_{0}(i)$ and $\r_{1}(i)\in S_{\n_{\a}\raj i}$.
By the definition of $R$, $\r_{2}(i)\in S_{\n_{\a}\raj i}$.
By ($\b$) above $\r_{2}(i)=\r_{1}(i)$ and $\r_{3}(i)=\r_{2}(i)$.
We can continue this and get $\r_{n}(i)=...=\r_{1}(i)$.
So $\n_{\a} (i)>\xi_{\b}(i)$.
$\eop$

\th 3.16 Lemma. There does not exist a sequence $(\t_{0},...\t_{n})$,
$n\in\o$, $n\ge 3$,
such that

(i) for all $m\le n$ there is active $\a$ such that $\t_{m}=\n_{\a}$
or $\t_{m}=\xi_{\a}$,

(ii) for all $m<n$ either

(a) $\t_{m}\overline{R}\t_{m+1}$

or

(b) there is
active $\a$ such that $\t_{m}=\n_{\a}$ and $\t_{m+1}=\xi_{\a}$
or $\t_{m}=\xi_{\a}$ and $\t_{m+1}=\n_{\a}$
\medskip
\noindent and at least case (b) exist in the sequence,

(iii) $\t_{0}=\t_{n}$,

(iv) for all
$m,m'<n$ if $m\ne m'$ then $\t_{m}\ne\t_{m'}$

\proof For a contradiction assume that such sequence exists.
By (ii) (b) we may choose the
sequence so that for some $\a$, $\t_{0}=\xi_{\a}$ and
$\t_{1}=\n_{\a}$. Then by (iv) and because
$n\ge 3$, $\t_{1}\overline{R}\t_{2}$.
By Lemma 3.12 $\n_{\a}\noR\xi_{\a}$ and so we may drop elements from
the sequence so that (i)-(iv) remain true, there are still at least
4 elements in the sequence and

(*)$\ \ \ \ $ if $m<n-1$ and $\t_{m}\overline{R}\t_{m+1}$ then
$\t_{m+1}\noR\t_{m+2}$.

By induction on $m<n$ we show that if $\t_{m}\noR\t_{m+1}$ then
$\t_{m+1}\overline{R}\t_{m+2}$ and if
$\t_{m}=\n_{\b}$ or $\t_{m}=\xi_{\b}$ for some $\b$ then
$\n_{\b}\nRm \xi_{\b}$.
Above we showed that $\n_{\a}\overline{R}\t_{2}$. By Lemma 3.14 (i)
$\t_{2}=\xi_{\b}$ for some active $\b$. By Lemma 3.14 (ii)
$\n_{\a}\nRm\xi_{\a}$. Then by (*) above $\t_{3}=\n_{\b}$.
By (iv) and Lemma 3.14 (i) $\t_{4}=\xi_{\g}$ for some active $\g$,
$\g\ne\b$ and $\n_{\b}\overline{R}\xi_{\g}$. By Lemma 3.14 (ii)
$\n_{\b}\nRm\xi_{\b}$. We can continue this and get the claim.

So there are active $\a_{0},...,\a_{m}$
such that the sequence is of the following form:
$$(\xi_{\a_{0}},\n_{\a_{0}},\xi_{\a_{1}},\n_{\a_{1}},...,
\n_{\a_{m}},\xi_{\a_{0}}).$$

We choose $i<\k$ so that for all $k\le m$, $i>j_{\a_{k}}$,
for all $k<m$, $i>l_{\a_{k}\a_{k+1}}$,
$i>l_{\a_{m}\a_{0}}$ and $i=\g +4p+3$ for some limit $\g$ and $p\in\o$.
By (vii')(b) $\xi_{\a_{o}}(i)\ge\n_{\a_{0}}(i)$. By
Lemma 3.15 $\n_{\a_{0}}(i)>\xi_{\a_{1}}(i)$. We can continue
this and finally we get $\n_{\a_{m}}(i)>\xi_{\a_{0}}(i)$. So
$\xi_{\a_{0}}(i)>\xi_{\a_{0}}(i)$, a contradiction. $\eop$

We define now $J_{\a}$, $\neg J_{\a}$, $K_{\a}$ and $\neg K_{\a}$
by induction on $\a <\a(*)$. We say that
$(J_{\a},\neg J_{\a},K_{\a},\neg K_{\a})$ is closed if

(i) $J_{\a}\cup K_{\a}$ and $\neg J_{\a}\cup\neg K_{\a}$ are
closed under $\overline{R}$,

(ii) if $\b$ is active then $\n_{\b}\in J_{\a}$ iff $\xi_{\b}\in\neg K_{\a}$
and $\n_{\b}\in\neg J_{\a}$ iff $\xi_{\b}\in K_{\a}$,

(iii) $J_{\a}\cap\neg J_{\a}=\empty$ and $K_{\a}\cap\neg K_{\a}=\empty$.

We assume that for all $\b <\a$ we have
defined $J_{\b}$, $\neg J_{\b}$, $K_{\b}$ and $\neg K_{\b}$
so that $(J_{\b},\neg J_{\b},K_{\b},\neg K_{\b})$ is closed.

If $\a$ is not active or
for some $\b <\a$, $\n_{\a}\in J_{\b}\cup\neg J_{\b}$
then we let $J_{\a}=\bigcup_{\b <\a}J_{\b}$,
$\neg J_{\a}=\bigcup_{\b <\a}\neg J_{\b}$,
$K_{\a}=\bigcup_{\b <\a}K_{\b}$ and
$\neg K_{\a}=\bigcup_{\b <\a}\neg K_{\b}$.

If $\a$ is active and
for all $\b <\a$, $\n_{\a}\not\in J_{\b}\cup\neg J_{\b}$
then we let $(J_{\a},\neg J_{\a},K_{\a},\neg K_{\a})$ be
such that it is closed and
$J_{\a}\supseteq\{\n_{\a}\}\cup\bigcup_{\b <\a}J_{\b}$,
$\neg J_{\a}\supseteq\bigcup_{\b <\a}\neg J_{\b}$,
$K_{\a}\supseteq\bigcup_{\b <\a}K_{\b}$ and
$\neg K_{\a}\supseteq\bigcup_{\b <\a}\neg K_{\b}$.
We prove the existence of these set by defining sets
$J^{i}_{\a}$, $\neg J^{i}_{\a}$, $K^{i}_{\a}$ and
$\neg K^{i}_{\a}$ by induction on $i<\vert\a (*)\vert^{+}$.

We let
$J^{0}_{\a}=\{\n_{\a}\}\cup\bigcup_{\b <\a}J_{\b}$,
$\neg J^{0}_{\a}=\bigcup_{\b <\a}\neg J_{\b}$,
$K^{0}_{\a}=\bigcup_{\b <\a}K_{\b}$ and
$\neg K^{0}_{\a}=\bigcup_{\b <\a}\neg K_{\b}$.
If $i<\vert\a (*)\vert^{+}$ is limit we let
$J^{i}_{\a}=\bigcup_{j <i}J^{j}_{\a}$ and similarly for the
other sets. If $i=j+1$ and odd then we let
the sets $J^{i}_{\b}$, $\neg J^{i}_{\b}$, $K^{i}_{\b}$ and $\neg K^{i}_{\b}$
be the least sets so that
$J^{i}_{\a}\supseteq J^{j}_{\a}$,
$\neg J^{i}_{\a}\supseteq\neg J^{j}_{\a}$,
$K^{i}_{\a}\supseteq K^{j}_{\a}$,
$\neg K^{i}_{\a}\supseteq\neg K^{j}_{\a}$
and $J^{i}_{\a}\cup K^{i}_{\a}$ and $\neg J^{i}_{\a}\cup\neg K^{i}_{\a}$ are
closed under $\overline{R}$.
If $i=j+1$ and even then if there is not active $\g$ such that

(1) $\n_{\g}\in J^{j}_{\a}$ and $\xi_{\g}\not\in\neg K^{j}_{\a}$ or

(2) $\n_{\g}\in\neg J^{j}_{\a}$ and $\xi_{\g}\not\in K^{j}_{\a}$ or

(3) $\xi_{\g}\in K^{j}_{\a}$ and $\n_{\g}\not\in\neg J^{j}_{\a}$ or

(4) $\xi_{\g}\in\neg K^{j}_{\a}$ and $\n_{\g}\not\in J^{j}_{\a}$

\noindent
then we let $J^{i}_{\a}=J^{j}_{\a}$ and similarly for the other sets.
Otherwise we let $\g$ be the least such ordinal and define

case (1): $J^{i}_{\a}=J^{j}_{\a}$, $\neg J^{i}_{\a}=\neg J^{j}_{\a}$,
$K^{i}_{\a}=K^{j}_{\a}$ and $\neg K^{j}_{\a}\cup\{\xi_{\g}\}$;

case (2): $J^{i}_{\a}=J^{j}_{\a}$, $\neg J^{i}_{\a}=\neg J^{j}_{\a}$,
$K^{i}_{\a}=K^{j}_{\a}\cup\{\xi_{\g}\}$ and $\neg K^{j}_{\a}$;

case (3): $J^{i}_{\a}=J^{j}_{\a}$, $\neg J^{i}_{\a}
=\neg J^{j}_{\a}\cup\{\n_{\g}\}$,
$K^{i}_{\a}=K^{j}_{\a}$ and $\neg K^{j}_{\a}$;

case (4): $J^{i}_{\a}=J^{j}_{\a}\cup\{\n_{\g}\}$,
$\neg J^{i}_{\a}=\neg J^{j}_{\a}$,
$K^{i}_{\a}=K^{j}_{\a}$ and $\neg K^{j}_{\a}$.

Finally we define $J_{\a}=\bigcup_{i<\vert\a (*)\vert^{+}}J^{i}_{\a}$
and similarly
for the other sets. If these sets are not as required then for some
$i=j+1<\vert\a (*)\vert^{+}$ even we have defined f.ex.
$\neg K^{i}_{\a}=\neg K^{j}_{\a}\cup\{\xi_{\g}\}$ while $\xi_{\g}$
belongs already to $K^{j}_{\a}$. If $i$ is the least such ordinal
then we can easily find a circle such that it contradicts Lemma 3.16.

So the sets $J_{\a}$, $\neg J_{\a}$, $K_{\a}$ and
$\neg K_{\a}$ exist.

We define $I_{0}=I^{-}_{0}\cup\bigcup_{\a <\a (*)}J_{\a}$ and
$I_{1}=I^{-}_{1}\cup\bigcup_{\a <\a (*)}K_{\a}$.

\th 3.17 Lemma. $I_{0}\not\cong I_{1}$.

\proof For a contradiction assume $g:I_{0}\rightarrow I_{1}$ is an
isomorphism. By Theorem 3.11 (iii) there exists an active $\a <\a (*)$
such that
for all $i\le\k$,
$$M^{\a}_{i}\prec (H_{<\k^{+}}(\l ),\in ,
I^{-}_{0},I^{-}_{1},g).$$

But then $\n_{\a}\in I_{0}$ iff $\xi_{\a}\not\in I_{1}$
and $g(\n_{\a})=\xi_{\a}$, which contradicts the assumption
that $g$ is an isomorphism. $\eop$

\th 3.18 Conclusion. Assume $\l=\mu^{+}$, $cf(\mu )=\mu$,
$\k=cf(\k )\le\mu$ and $\l^{<\k}=\l$. Then there are $\l^{+},\k +1$-trees
$I_{0}$ and $I_{1}$
such that $I_{0}\not\cong I_{1}$ and
$$I_{0}\equiv^{\l}_{\mu\times\k}I_{1}.$$
If $\l^{\k}=\l$ then $I_{0}$ and $I_{1}$ are of cardinality $\l$.

Notice that if we replace Theorem 3.11 with a slightly stronger
black box (see [Sh3]), we can, instead of two $\l^{+},\k$-trees,
get $2^{\l}$ $\l^{+},\k$-trees such that any two of them
satisfy Conclusion 3.18.

\chapter{4. On structure of trees of fixed height}

In this chapter we will show that trees of fixed height
are isomorphic
if they are equivalent up to some relatively small tree. This implies
that essentially the same is true for the models of the canonical example of
unsuperstable theories (see [HT]).

\th 4.1 Definition. ([Sh1])
Let $\l$ be a regular cardinal. We define $I[\l ]$
to be the set of $A\subseteq\l$ such that there exist a cub $E\subseteq\l$
and $\P =\{ P_{\a}\vert\ \a<\l\}$ satisfying

(i) $P_{\a}$ is a set of
subsets of $\a$ and $\vert P_{\a}\vert <\l$;

(ii) for all limit $\d\in A\cap E$ such that $cf(\d )<\d$, there exists
$C\subseteq\d$
such that

(a) the order type of $C$ is $<\d$ and $sup\ C=\d$;

(b) $C\cap\a\in\bigcup_{\b <\d}P_{\b}$ for all $\a <\d$.

Notice that for example $\o_{1}\in I[\o_{1}]$: Let $E\subset\o_{1}$
be the set of all limit ordinals $<\o_{1}$ and $\P =\{ P_{\a}\vert\ \a<\l\}$
such that $P_{\a}=\{ B\subseteq\a\vert\ \vert B\vert <\o\}$. Then
(i) and (ii) above are satisfied.
For further properties of $I[\l ]$ see [Sh1].

\th 4.2 Definition. Let $\l$ be a regular cardinal and $t$ a $\l^{+},\l$-tree
of cardinality $\l$. Let $\{ x_{i}\vert\ i<\l\}$ be an enumeration of
$t$ and let $t'$ be a subtree of $t$.
Then $S[t']$ is the set of those limit ordinals $\d <\l$ which
satisfy the following condition (*):

(*) $\{ x_{i}\in t'\vert\ i<\d\}$ contains a branch of length $\d$.

\relax From now on we assume that when ever we talk about a tree $t$, we have
fixed an enumeration $\{ x_{i}\vert\ i<\vert t\vert\}$ for it. We
assume that the enumeration is such that if $x_{i}<x_{j}$ then
$i<j$.

\th 4.3 Definition. Let $\l$ and $\k$ be regular cardinals, $\k <\l$
and $t$ a $\l^{+},\l$-tree of cardinality $\l$.
Let $\{ x_{i}\vert\ i<\l\}$ be the enumeration of
$t$.
We say that $t$ is
$\l ,\k$-large if $t$ satisfies the following condition: There are
sets $E_{\xi}$, $\xi\le\k$, such that

(i) $E_{\xi}\subseteq t$ and if $\xi\ne\xi'$ then
$E_{\xi}\cap E_{\xi'}=\empty$;

(ii) for $\xi <\d$ and $x\in E_{\d}$ there is a unique $y\in E_{\xi}$
such that $y<x$;

(iii) if $\d\le\k$ is limit, $x_{\xi}\in E_{\xi}$ for all $\xi <\d$ and
$(x_{\xi})_{\xi <\d}$ is increasing then there is $y\in E_{\d}$ such
that $x_{\xi}<y$ for all $\xi <\d$;

(iv) if $\xi <\k$, $x\in E_{\xi}$ then we write
$$t_{x}=\{ y\in t\vert\ x\le y\ \hbox{\sl and there is}\ z\in E_{\xi +1}
\ \hbox{\sl such that}\ y< z\}$$
and require
than there exists a set $\T$ of regular cardinals $<\l$ such that

(a) $S[t_{x}]\cup\{\d<\l\vert\ cf(\d )<\d,\ cf(\d )\in\T\}$ contains a
cub set (in $\l$);

(b) $\{\d <\l\vert\ cf(\d )<\d ,\ cf(\d )\in\T,\ \d\not\in
S[t_{x}]\}\in I[\l ]$;

(c) for $\d\in\T$ there is $y\in t_{x}$ such that
the order type of $\{ z\vert\ x\le z< y\}$ is $\d$;

(v) if $\g =\b +1<\k$,
$(x_{\xi})_{\xi <\d}$ is an increasing sequence in $t$,
$x_{0}\in E_{\b}$ and for
all $\xi <\d$ there is $y_{\xi}\in E_{\g}$ such that $x_{\xi}<y_{\xi}$, then
there is $y\in E_{\g}$ such that $x_{\xi}<y$ for all $\xi <\d$.

Notice that if $\l =\mu^{+}$, $\l\in I[\l]$
and $\k <\l$ is regular then $\mu\times\k +1$
is a $\l ,\k$-large $\l^{+},\l$-tree.
If $\l$ is weakly compact then
there is no $\l ,\k$-large $\l^{+},\l$-trees.

The proof of the theorem below is a modification of the proof
of related result in [HT]. The most conspicuous difference is
the use of elementary submodels of $H(\l^{*})$. They are used
only to make it easier to define the closures needed in the
proof.

\th 4.4 Theorem. Let $\l$ and $\k$ be regular cardinals, $\k <\l$ and
$I_{0}$ and $I_{1}$ be $\l^{+},\k +1$-trees.
Assume $t$ is a $\l,\k$-large
$\l^{+},\l$-tree of cardinality $\l$. Then
$$I_{0}\equiv^{\l}_{t}I_{1}\ \ \Leftrightarrow\ \ I_{0}\cong I_{1}.$$

\proof Without loss of generality we may assume that $I_{0}$ and $I_{1}$
are such that if $x,y\in I_{0}$ ($\in I_{1}$), they have no immediate
predecessors, $x\sim y$ and $pred(x)$ is of power $<\k$ then $x=y$.

Let $\r$ be a winning strategy of $\exists$ in
$G^{\l}_{t}(I_{0},I_{1})$.
We define by induction on $\a \le\k$ the following:

(i) an isomorphism $f_{\a}$ from
$I^{\le\a}_{0}$
onto $I^{\le\a}_{1}$;

(ii) for each $x\in I^{\le\a}_{0}\cup I^{\le\a}_{1}$ we define
an initial segment $R_{x}=((a_{i},X_{i},Y_{i}))_{i\le\b}$ of a
play in $G^{\l}_{t}(I_{0},I_{1})$,
such that $x\in\bigcup_{i\le\b}(rng(X_{i})\cup rng(Y_{i}))$,
$rng(X_{i})\cup rng(Y_{i})\subseteq I^{\le\a}_{0}\cup I^{\le\a}_{1}$
for all $i<\b$, $\exists$ has used
$\r$ and if $x$ is not a leaf then for some $\d <\k$
there is $a_{x}\in E_{\d}$ such that
$a_{i}\le a_{x}$ for all $i<\b$.
Furthermore we require that
if $x<x'$ then $R_{x}$ is an initial
segment of $R_{x'}$ and for each $x\in I^{\le\a}_{0}$
$f_{\a}(x)$ is the element $\exists$ has
chosen to be the image of $x$ in $R_{x}$.

If we can do this we have clearly proved the theorem. The cases
$\a=0$ and $\a$ is limit are trivial. So we assume that
$\a=\g +1$.

Let $z\in I^{\le\g}_{0}-\cup_{\d <\g}I^{\le\d}_{0}$. Clearly it
is enough to define $f_{\a}\raj succ(z)$ and $R_{x}$ for all
$x\in succ(z)$ so that $f_{\a}\raj succ(z)$ is onto
$succ(f_{\g}(z))$.
Let $y=f_{\g}(z)$ and let $n:\l\rightarrow t$
be the function that gives the enumeration of $t$, $t=\{ n(i)\vert\ i<\l\}$
(see the assumption after Definition 4.2).
Let $R_{z}=((a_{i},X_{i},Y_{i}))_{i\le\b}$. By induction assumption
there is
$a_{z}\in E_{\d}$, $\d <\k$, such that
$a_{i}< a_{z}$ for all $i\le\b$.
Let $E$ and
$\P =\{ P_{i}\vert\ i<\l\}$ be the sets which show that
$$\{\d <\l\vert\ cf(\d )<\d , cf(\d )\in\T,\ \d\not\in
S[t_{a_{z}}]\}\in I[\l ].$$

Let $\l^{*}$ be large enough, say $(\beth_{10}(\l))^{+}$. We choose
$\A_{i}$, $i<\l$, so that

(a) $\vert\A_{i}\vert <\l$ and
$\A_{i}\prec (H(\l^{*}),\in ,I_{0},I_{1},t,<_{0},<_{1},<)$, where
$<_{0}$ denotes the ordering of $I_{0}$, $<_{1}$ denotes the
ordering of $I_{1}$ and $<$
denotes the ordering of $t$;

(b) $\r ,n,(E_{\xi}\vert\ \xi\le\k ),
E,(P_{i}\vert\ i<\l ),R_{z},\l ,\b ,a_{z}\in\A_{0}$,
$\k +1\subseteq\A_{0}$ and $i\subseteq\A_{i}$;

(c) $\A_{i}\prec\A_{j}$ if $i<j$ and $\A_{i}=\cup_{j<i}\A_{j}$ if
$i$ limit;

(d)
for all $i\le\b$, $dom(X_{i})\in\A_{0}$ (see Definition
2.2);

(e) $\A_{i}\cap\l$ is ordinal, $\A_{i}\in\A_{i+1}$ and
$\A_{i}\cap\l\in\A_{i+1}$;

(f) $succ(z)\cup succ(y)\subseteq\bigcup_{i<\l}\A_{i}$;

(g) if $x\in t\cap\A_{i}$, $y\in t$ and $y<x$, then $y\in\A_{i}$.

Let
$$C\subseteq S[t_{a_{z}}]\cup\{\d <\l\vert\ cf(\d )<\d
\ \hbox{\rm and}\ cf(\d )\in\T\}$$
be cub.
We may assume that for all $c\in C$, $\A_{c}\cap\l =c$ and $c\in E$.

For all $i<\l$ we define by induction $c_{i}\in C$ and
$f_{\a}\raj (succ(z)\cap \A_{c_{i}})$. If $i$ is limit
then $c_{i}=\bigcup_{j<i}c_{j}$ and
$f_{\a}\raj (succ(z)\cap \A_{c_{i}})$ is already defined.

Assume that we have defined $c_{i}$ and
$f_{\a}\raj (succ(z)\cap \A_{c_{i}})$
as wanted and
$$rng(f_{\a}\raj (succ(z)\cap \A_{c_{i}}))=succ(y)\cap\A_{c_{i}}.$$
Let us define $c_{i+1}$ and
$$f_{\a}\raj (succ(z)\cap (\A_{c_{i+1}}-\A_{c_{i}})).$$

Now either $c_{i}\in S[t_{a_{z}}]$ or $c_{i}\in\{\d <\l\vert\ cf(\d )<\d
\ \hbox{\rm and}\ cf(\d )\in\T\}$.

(1) $c_{i}\in S[t_{a_{z}}]$: Let $B\in\A_{c_{i}+1}$ be a branch in
$$S[t_{a_{z}}]\cap\A_{c_{i}}=\{ n(j)\vert\ j<c_{i}\}$$
of length $c_{i}$.
Let $h\in\A_{c_{i}}$ be a one-one
function from $(succ(z)\cup succ(y))\cap\A_{c_{i}}$ to
$\A_{c_{i}}\cap\l$.
We let the players continue the play
$R_{z}$ so that in the next $c_{i}$ moves $\forall$ chooses the sets
$\{ h^{-1}(\d )\}$, $\d <c_{i}$, from $I_{0}\cup I_{1}$ and from
$t$ he chooses elements of $B$. We let $\exists$ follow $\r$.
If $B'$ is an initial segment of $B$ then
$B'=\{ y\in t\vert\ a_{z}\le y<x\}$ for some $x\in B$. So $B'\in\A_{c_{i}}$,
which implies that every initial segment of the play belongs to $\A_{c_{i}}$.
Because $\A_{c_{i}}$ is closed under $\r$,
all the elements $\exists$ chooses are from $\A_{c_{i}}$. It is also
easy to see that this play belongs to $\A_{\g}$ for all $\g >c_{i}$.

By Definition 4.3 (v) we can find $a\in E_{\d +1}\cap\A_{c_{i}+1}$,
such that $a$ is larger than any element $b\in t$
chosen by $\forall$ in the play above. Let
$$C'\subseteq S[t_{a}]\cup\{\d <\l\vert\ cf(\d )<\d
\ \hbox{\rm and}\ cf(\d )\in\T\}$$
be cub. Let $c_{i+1}\in C\cap C'$ be such that $c_{i+1}>c_{i}$.
Then $a\in\A_{c_{i+1}}$.
Now either $c_{i+1}\in S[t_{a}]$ or $c_{i+1}\in\{\d <\l\vert\ cf(\d )<\d
\ \hbox{\rm and}\ cf(\d )\in\T\}$. In the first case we let
$\forall$ play the elements $(succ(z)\cup succ(y))\cap\A_{c_{i+1}}$
as above. So let us assume that $c_{i+1}\not\in S[t_{a}]$ and
$c_{i+1}\in\{\d <\l\vert\ cf(\d )<\d
\ \hbox{\rm and}\ cf(\d )\in\T\}$. Especially then
$$(*)\ \ \ \ c_{i+1}\in E\cap \{\d <\l\vert\ cf(\d )<\d ,\ cf(\d )
\in\T,\ \d\not\in
S[t_{a}]\} .$$

Let $h'\in\A_{c_{i+1}}$ be a one-one
function from $(succ(z)\cup succ(y))\cap\A_{c_{i+1}}$ to
$c_{i+1}=\A_{c_{i+1}}\cap\l$. Let
$$D'\subseteq c_{i+1}$$
be a set such that for all $\xi <c_{i+1}$, $\xi\cap D'\in
\bigcup_{j<c_{i+1}}P_{j}$,
$sup\ D'=c_{i+1}$ and the order type of $D'$ is $cf(c_{i+1})$. The
existence of this set follows from (*) above.
Let $D=\{ d_{j}\vert\ j<cf(c_{i+1})\}$ be the closure of $D'$ in $c_{i+1}$.
Because $cf(c_{i+1})\in\A_{c_{i+1}}$, it is easy to see that in
$t_{a}\cap\A_{c_{i+1}}$ there is a branch $B$ of length $cf(c_{i+1})$.
We let the players continue the play above
so that in the next $cf(c_{i+1})$ moves $\forall$ chooses the sets
$\{ h'^{-1}(k)\vert\ k<d_{j}\}$ from $I_{0}\cup I_{1}$,
$j<cf(c_{i+1})$,
and from $t$ he chooses elements of $B$.
We let $\exists$ follow $\r$.

Because $\bigcup_{i<c_{i+1}}P_{i}\subseteq\A_{c_{i+1}}$, every initial
segment of this play is in $\A_{c_{i+1}}$ and so all elements chosen
by $\exists$ from $I_{0}\cup I_{1}$
are from $\A_{c_{i+1}}$. Then
by using the moves of
$\exists$ we can define
$$f_{\a}\raj (succ(z)\cap (\A_{c_{i+1}}-\A_{c_{i}})).$$
For each $x\in succ(z)\cap (\A_{c_{i+1}}-\A_{c_{i}})$,
$R_{x}$ will be the play defined above.

(2) $c_{i}\not\in S[t_{a_{z}}]$:
Now we first let $\forall$ play
the elements of $(succ(z)\cup succ(y))\cap\A_{c_{i}}$ as in the
second half of the case (1) and then continue as above.
Notice that in this case (also) we have to define the first
$cf(c_{i})$ moves so that the play belongs to $\A_{c_{i+1}}$.
We can guarantee this by choosing $D'\subseteq c_{i}$ so
that $D'\in\A_{c_{i}+1}$.
$\eop$

\th 4.5 Remark. Let $\l =\mu^{+}$ and $\k <\l$ regular. Let
$I_{0}$ and $I_{1}$ be $\l^{+},\k+1$-trees. Assume $\l\in I[\l ]$.
Above we proved that
if $\a =\mu\times\k +1$ then
$$(*)\ \ \ \ I_{0}\equiv^{\l}_{\a}I_{1}\ \ \Leftrightarrow
\ \ I_{0}\cong I_{1}.$$
In Chapter 3 we showed that if $\mu$ is regular then this is
best possible. But if $\mu$ is not regular then we can get better
results.

If $\k <cf(\mu )<\mu$ then (*) is true if $\a =\mu$ and
if $\k =cf(\mu )<\mu$ then (*) is true if $\a =\mu +1$. This
can be proved as Theorem 4.4.

\vfill
\eject

\chapter{References}

\item{[HS]} T.Hyttinen and S.Shelah, Constructing strongly equivalent
nonisomorphic models for unsuperstable theories, part B, to appear.

\item{[HT]} T.Hyttinen and H.Tuuri, Constructing strongly equivalent
nonisomorphic models for unstable theories, APAL 52, 1991, 203--248.

\item{[Je]} T.Jech, Set Theory, Academic Press, New York, 1978.

\item{[MSV]} A. Mekler, S. Shelah and J. V""n"nen, The Ehrenfeucht-Fraisse-game
of length $\o_{1}$, to appear.

\item{[Sh1]} S.Shelah, Nonelementary classes: Part II, in: J.T.Baldwin (ed.),
Classification Theory, Springer Lecture Notes 1292, (Springer-Verlag, 1987)
419--497.

\item{[Sh2]} S.Shelah, Classification Theory, Stud. Logic Found. Math.
92 (North-Holland, Amsterdam, 2nd rev. ed., 1990).

\item{[Sh3]} S.Shelah, Non-structure Theory, to appear.
\bigskip

Tapani Hyttinen

Department of Mathematics

P. O. Box 4

00014 University of Helsinki

Finland
\bigskip

Saharon Shelah

Institute of Mathematics

The Hebrew University

Jerusalem

Israel
\medskip

Rutgers University

Hill Ctr-Busch

New Brunswick

NJ 08903

USA

\end